\documentclass[a4paper,12pt]{article}
\usepackage{amsmath}
\usepackage{amssymb,eurosym}
\usepackage{mathtools}
\usepackage[dvips]{graphicx}
\usepackage{lscape}
\usepackage{color}
\usepackage{graphics}
\usepackage{graphicx}
\graphicspath{{./Graphiken/}}
\usepackage{amsmath}
\usepackage[misc]{ifsym}
\usepackage[nospace,noadjust]{cite}
\usepackage[pdftoolbar=true,
            pdfmenubar=true,
            pdfpagemode=UseOutlines,
            bookmarksnumbered=true,
            linktocpage=true,
            colorlinks=false]{hyperref}
\usepackage{verbatim}
\newtheorem{definition}{Definition}[section]
\newtheorem{theorem}[definition]{Theorem}
\newtheorem{lemma}[definition]{Lemma}
\newtheorem{corollary}[definition]{Corollary}
\newtheorem{proposition}[definition]{Proposition}
\newtheorem{remark}[definition]{Remark}
\newtheorem{example}[definition]{Example}

\newtheorem{assumption}[definition]{Assumption}

\newtheorem{conjecture}[definition]{Conjecture}

\newcommand{\bdefi}{\begin{definition}}
\newcommand{\edefi}{\end{definition}}
\DeclareMathOperator{\gem}{Gem}
\DeclareMathOperator{\new}{New}
\newcommand{\blem}{\begin{lemma}}
\newcommand{\elem}{\end{lemma}}
\newcommand{\bthe}{\begin{theorem}}
\newcommand{\ethe}{\end{theorem}}
\newcommand{\bcor}{\begin{corollary}}
\newcommand{\ecor}{\end{corollary}}
\newcommand{\bprop}{\begin{proposition}}
\newcommand{\eprop}{\end{proposition}}
\newcommand{\brem}{\begin{remark}}
\newcommand{\erem}{\end{remark}}
\newcommand{\bex}{\begin{example}}
\newcommand{\eex}{\end{example}}
\newcommand{\bass}{\begin{assumption}}
\newcommand{\eass}{\end{assumption}}
\newcommand{\bconj}{\begin{conjecture}}
\newcommand{\econj}{\end{conjecture}}

\renewcommand{\i}{\begin{itemize}}
\newcommand{\ii}{\end{itemize}}
\newcommand{\q}{\begin{equation}}
\newcommand{\qq}{\end{equation}}
\newcommand{\qa}{\begin{eqnarray*}}
\newcommand{\qqa}{\end{eqnarray*}}
\renewcommand{\a}{\begin{array}}
\renewcommand{\aa}{\end{array}}

\newcommand{\1}{{1\!\!1}}

\DeclareMathOperator{\conv}{conv}

\renewcommand{\epsilon}{\varepsilon}

\renewcommand{\L}{{\cal L}}

\newcommand{\N}{\mathbb N}
\newcommand{\C}{\mathbb C}

\newcommand{\G}{{\cal G}}
\renewcommand{\phi}{\varphi}
\newcommand{\pmat}[1]{\begin{pmatrix}#1\end{pmatrix}}

\newcommand{\proof}{{\bf Proof}. }
\newcommand{\qed}{~\hfill$\bullet$}
\newcommand{\R}{\mathbb R}

\DeclareMathOperator{\sign}{sign}

\renewcommand{\subset}{\subseteq}

\DeclareMathOperator{\vertex}{vert}

\numberwithin{equation}{section}

\begin{document}

\title{\bf On Globally Diffeomorphic Polynomial Maps via Newton Polytopes and Circuit Numbers}

\author{
Tom\'a\v{s} Bajbar\thanks{Institute of Operations Research, Karlsruhe Institute of Technology (KIT), Germany, bajbar@kit.edu (\Letter)}
\and
Oliver Stein\thanks{Institute of Operations Research, Karlsruhe Institute of Technology (KIT), Germany, stein@kit.edu}
} 


\maketitle

\begin{abstract} In this article we analyze the global diffeomorphism property of polynomial maps $F:\R^n\to\R^n$ by studying the properties of the Newton polytopes at infinity corresponding to the sum of squares polynomials $\|F\|_2^2$. This allows us to identify a class of polynomial maps $F$ for which their global diffeomorphism property on $\R^n$ is equivalent to their Jacobian determinant $\det JF$ vanishing nowhere on $\R^n$. In other words, we identify a class of polynomial maps for which the Real Jacobian Conjecture, which was proven to be false in general, still holds. 
\end{abstract}

\bigskip
\bigskip
\bigskip
\bigskip
\bigskip



\newpage

\section{Introduction}\label{sec:intro} 

It is an interesting question how to verify or disprove whether a given differentiable map $F:\R^n\to\R^n$ is globally invertible with a differentiable inverse $F^{-1}:\R^n\to\R^n$. In the present work, we shall call such maps global diffeomorphisms of $\R^n$ onto itself. The first well-known characterization of this global diffeomorphism property dates back to the work of Hadamard \cite{Had1,Had2,Had3} and states that it is equivalent to the determinant $\det JF$ of the Jacobian matrix $JF$ of $F$ vanishing nowhere on $\R^n$, and to $F$ being proper (cf. Th.~\ref{thm:Had} below). Here, $F$ is called proper if preimages of compact sets under $F$ always are compact. 

In the case of complex polynomial maps $F:\C^n\to\C^n$, the characterization of their global invertibility property directly refers to the Jacobian Conjecture from algebraic geometry, first formulated in \cite{Kel} and asserting that if $\det JF$ is a nonzero constant function then $F$ possesses a global polynomial inverse $F^{-1}:\C^n\to\C^n$. There exists a vast number of partial results on this conjecture where different approaches are used (see, e.g. \cite{Wa,BCR,BR,Y, Rab, Moh}). For more details on this open problem we refer to the survey papers \cite{Dr, Wr, Ess1,Ess2}. 

Following \cite{Kur}, in the setting of real polynomial maps $F:\R^n\to\R^n$, the injectivity of $F$ implies its surjectivity \cite{BR}, and the global inverse $F^{-1}$ of $F$ is a polynomial if and only if $\det JF$ is a nonzero constant function \cite{BCR}. If merely the existence, but not necessarily the polynomiality of the inverse map $F^{-1}$ is sought for, one may conjecture that $\det JF$ vanishing nowhere on $\R^n$ implies the injectivity of $F$, and hence, also the existence of its global inverse $F^{-1}:\R^n\to\R^n$. This is the so-called ``Real Jacobian Conjecture''. It was, however, proven to be false by Pinchuk in \cite{Pi}, where a counterexample of a non-injective polynomial map $F:\R^2\to\R^2$ is constructed with  $\det JF$ vanishing nowhere on $\R^2$. 

Since, by Hadamard's above-mentioned theorem, the non-vanishing property of $\det JF$ actually is necessary for the global diffeomorphism property of $F$, it is thus an interesting question which additional conditions imposed on $F$, general enough, can assure that $F$ is a global diffeomorphism of $\R^n$ onto itself. Answering this question, which is also posed by Bivi\`{a}-Ausina in \cite{Biv} and which is of significant importance in \cite{Tib} as well, is the main motivation for the present article. 

From Hadamard's theorem it is clear that such additional conditions must be related to the properness of $F$. Since the latter properness may be characterized by the coercivity of the sum of squares polynomial $\|F\|_2^2$, sufficient conditions for the coercivity of polynomials will be the main tools used in the present article. In fact, it will turn out that for a broad class of polynomial maps $F:\R^n\to\R^n$ the
coercivity of $\|F\|_2^2$ follows from $\det JF$ being nonzero on $\R^n$, so that at least for this class of polynomial maps the Real Jacobian Conjecture turns out to be true.

We mention that, after Hadamard's fundamental contribution, further important results on the global diffeomorphism property were proved by Levy \cite{Levy}, Banach and Mazur \cite{BM}, Caccioppoli \cite{Cac}, Plastock \cite{Plas} and Rabier \cite{Rab2}. For a brief summary and further details see, e.g., \cite{Gor, Rab2, Zam}. It is worth mentioning that in mathematical economics, global invertibility properties of maps, as an object of interest, was originally highlighted in \cite{Sam} and subsequently studied in \cite{GaNi, Mas, Chich}. The global invertibility of homogeneous maps, which are not smooth at the origin, is studied in \cite{Ruz}.

This article is structured as follows. In Section \ref{sec:coerc} the notation and the main results from \cite{BS} on coercive polynomials and their Newton polytopes are briefly summarized. 
More precisely, we will first recall the definition of the gem of a multivariate polynomial $f\in\R[x]$, a geometric structure associated to the Newton polytope at infinity $\new_{\infty}(f)$, which enables one to define the broad class of so-called gem regular polynomials (see Def.~\ref{def:Fregular} below). Next we recall the main theorem from \cite{BS}, which characterizes the coercivity of  gem regular polynomials $f\in\R[x]$ via three conditions \eqref{eq:c1}-\eqref{eq:c3} imposed on the vertex set at infinity of the corresponding Newton polytopes at infinity $\new_{\infty}(f)$ (see Th.~\ref{the:char} below). For polynomials $f\in\R[x]$ which are not necessarily gem regular, we further recall some important necessary and some sufficient conditions for coercivity (see Ths.~\ref{the:necFirreg} and \ref{the:suffdeg} below), where so-called circuit numbers appear. 

In Section \ref{sec:GlobDiff}, we show that every sum of squares polynomial $\|F\|_2^2$ corresponding to some polynomial map $F:\R^n\to\R^n$ fulfills conditions \eqref{eq:c1} and \eqref{eq:c2} and, using a determinant formula for Jacobians $JF$ (see Lem.~\ref{thm:det} below), we prove that polynomials $\|F\|_2^2$ corresponding to polynomial maps $F$ with nonvanishing Jacobian determinants $\det JF$ fulfill also the condition \eqref{eq:c3} (see Props.~\ref{lem:main} and \ref{prop:main} below). Finally, a combination of Hadamard's theorem (see Th.~\ref{thm:Had} below) and the coercivity results from Section \ref{sec:coerc} enables us to identify a class of polynomial maps $F:\R^n\to\R^n$ whose global diffeomorphism property on $\R^n$ is equivalent to their Jacobian determinant $\det JF$ vanishing nowhere on $\R^n$, which is the main result of the present paper (see Ths.~\ref{thm:main} and \ref{thm:maindeg} below). 

This class of polynomial maps $F:\R^n\to\R^n$ is described in terms of so-called Newton polytopes at infinity $\new_{\infty}(\|F\|_2^2)$ corresponding to $\|F\|_2^2$. More precisely, for a given polynomial map $F:\R^n\to\R^n$, in order to verify whether $F$ belongs to the latter class, one has to identify the vertex set at infinity $V(\|F\|_2^2)$, the set of so-called gem-degenerate exponent vectors $D(\|F\|_2^2)$, and for the latter one also has to compute the corresponding circuit numbers (for definitions, see Section \ref{sec:coerc}). The first may be realized by, for example, vertex- or facet enumeration algorithms (for more details see, e.g., \cite{AF, BFM}). 

We illustrate our main results in Example \ref{ex:thm2}, where a one-parametric family of polynomial diffeomorphisms of $\R^2$ onto itself is analyzed by using our techniques. Since for some singular parameter value
these techniques are not directly applicable, in Section~\ref{sec:lin} we also prove the invariance of the coercivity property under linear coordinate transformations, 
and show that our main results may be generalized by replacing the assumptions on $\|F\|_2^2$ by assumptions on $\|F\circ A^{-1}\|_2^2$ for some regular matrix $A\in\R^{n\times n}$ (see Cors.~\ref{cor:main} and \ref{cor:maindeg} below). In Example \ref{ex:thm2_voll}, we 
use such a transformation to apply our techniques to treat the case of the singular parameter from Example \ref{ex:thm2}.
The article closes with some final remarks in Section~\ref{sec:fin}.

\section{Review of results on coercive polynomials and their Newton polytopes at infinity}\label{sec:coerc}

First, we introduce some notation we shall use throughout the present article. Let $\R[x] =\R[x_1,\dots,x_n]$ denote the ring of polynomials  in $n$ variables with real coefficients. For $f\in\R[x]$ we write $f(x)=\sum_{\alpha\in A(f)}f_{\alpha}x^\alpha$ with $A(f)\subseteq\N_0^n$, $f_{\alpha}\in\R\setminus\{0\}$ for $\alpha\in A(f)$, and $x^{\alpha}=x_1^{\alpha_1}\dots x_n^{\alpha_n}$ for $\alpha\in\N_0^n$. 

A function $f:\R^n\rightarrow\R$ is called coercive on $\R^n$, if $f(x)\rightarrow +\infty$ holds whenever $\|x\|\rightarrow +\infty$, where $\| \cdot \|$ denotes some norm on $\R^n$. In \cite{BS} it is shown how the coercivity of multivariate polynomials can often be analyzed by studying so-called Newton polytopes at infinity, whose definition we recall in the next step.

For $f\in\R[x]$, the set $\new_{\infty}(f):=\conv \left(A(f)\cup\{0\}\right)$, that is, the convex hull of the set $A_0(f):=A(f)\cup\{0\}$ is called the \textit{Newton polytope at infinity} of the polynomial $f$, and the set $\new (f):=\conv \left(A(f)\right)$ is called the \textit{Newton polytope} of the polynomial $f$. The set $V_0(f)$ denotes the set of all vertices of $\new_{\infty}(f)$. With $\mathbb{H}:=\{h\in\R|\,h\geq 0\}$ one obtains, due to $0\in\new_{\infty}(f)\subseteq\mathbb{H}^n$, that the set $V_0(f)$ always contains the origin. For the later purposes of this work we shall define the vertex set of $\new_{\infty}(f)$ at infinity $V(f):=V_0(f)\setminus\{0\}$  and the set $V^c(f):=A(f)\setminus V(f)$ of all exponent vectors of $f$ which are no vertices at infinity of $\new_{\infty}(f)$.

Various algebraic and analytic properties of polynomials are encoded in the properties of their Newton polytopes. To name some of them, for example the number of isolated roots of $n$ polynomial equations in $n$ unknowns can be bounded by the (mixed) volumes of their Newton polytopes (cf., e.g., \cite{KaKu,Kus,Sturmfels}), absolute irreducibility of a polynomial is implied by the indecomposability of its Newton polytope in the sense of Minkowski sums of polytopes \cite{SGao}, and there are also some results dealing with Newton polytopes in elimination theory \cite{Elimin}. 

In the following the main concepts and results on coercive polynomials and their Newton polytopes at infinity from \cite{BS} are briefly recalled.

For a polynomial $f\in\R[x]$ let
\[
\mathcal{G}(f)\ :=\ \{G\subset\R^n|\ G\neq\emptyset\ \text{is a face of}\ \new_{\infty}(f)\ \text{with}\ 0\not\in G\}
\]
be the set of all nonempty faces of  $\new_{\infty}(f)$ not including the origin. The set
\[
\gem(f)\ :=\ \bigcup_{G\in\mathcal{G}(f)}\,G,
\] 
is called the \emph{gem} of $f$ (in \cite{Kou, Pha08, Tha13} also called the ``Newton boundary at infinity'') and gives rise to the following important regularity concept for polynomials:


\begin{definition}[\mbox{\cite{BS}}]\label{def:Fregular}  Let $f\in\R[x]$ be given.
\begin{itemize}
\item[a)] An exponent vector $\alpha\in A(f)$ is called \emph{gem degenerate} if  $\alpha\in V^c(f)\cap G$ holds for some $G\in\mathcal{G}(f)$. We denote the set of all gem degenerate points $\alpha\in A(f)$ by $D(f)$.
\item[b)] The polynomial $f$ is called \emph{gem re\-gular} if the set $D(f)$ is empty,
otherwise it is called \emph{gem irregular}. 
\end{itemize}
\end{definition}

Clearly, gem regularity of $f\in\R[x]$ is equivalent to $V^c(f)\cap G=\emptyset$ for all $G\in\cal{G}$.
Furthermore, the definition of $D(f)$ gives rise to a partitioning of $V^c(f)$ into $D(f)$ and a set of ``remaining exponents''
$R(f):=V^c(f)\setminus D(f)$, so that we may write
\q\label{eq:VDR}
A(f)\ =\ V(f)\ \dot\cup\ D(f)\ \dot\cup\ R(f).
\qq
Using \eqref{eq:VDR} together with the notation $f^W(x):=\sum_{\alpha\in W}f_\alpha x^\alpha$ for some $W\subseteq A(f)$, any $f\in\R[x]$ can be expressed as
\q\label{eq:fVDR}
f=f^{V(f)}+f^{D(f)}+f^{R(f)}.
\qq

The following three conditions from \cite{BS} are crucial for analyzing the coercivity of $f\in\R[x]$ on $\R^n$. Here and subsequently we put $I:=\{1,\dots,n\}$.
\q
V(f)\subset 2\N_0^n.\tag{C1}\label{eq:c1}
\qq
\q
\text{All $\alpha\in V(f)$ satisfy $f_\alpha>0$.}\tag{C2}\label{eq:c2}
\qq
\q
\text{For all $i\in I$ the set $V(f)$ contains a vector of the form $2k_ie_i$ with $k_i\in\N$.}\tag{C3}\label{eq:c3}
\qq

\bthe[Characterization of coercivity, \mbox{\cite[Th. 3.2]{BS}}]\label{the:char}\ \\
Let $f\in\R[x]$ be gem regular. Then the following assertions are equivalent.

\i
\item[a)] $f$ is coercive on $\R^n$.
\item[b)] $f$ fulfills the conditions \eqref{eq:c1}-\eqref{eq:c3}.
\ii
\ethe

To treat the gem irregular case, recall that, by Carath\'eodory's theorem, for any degenerate exponent vector $\alpha^\star\in D(f)$ there exists a set of affinely independent points
$V^\star\subset V(f)$ with $\alpha^\star\in\conv V^\star$. In the case that a simplicial face $G\subset\mathcal{G}(f)$ contains $\alpha^\star$, the
set $V^\star$ can be chosen as the vertex set $V_G$ of $G$. For non-simplicial faces $G$, however, there may exist several possibilities
to choose $V^\star\subset V_G$.

For any set of affinely independent points $V^\star$ with $\alpha^\star\in\conv V^\star$, the solution $\lambda$ of
\[
\sum_{\alpha\in V^\star}\lambda_\alpha\pmat{\alpha\\1}\ =\ \pmat{\alpha^\star\\1},\quad \lambda_\alpha\ge0,\ \alpha\in V^\star
\]
is unique
and, using the natural convention $0^0:=1$ in the polynomial setting (to cover the case of vanishing coefficients $\lambda_\alpha$), we may define the \emph{circuit number} 
(cf. \cite{Ili14Wa}) 
\[
\Theta(f,V^\star,\alpha^\star)\ =\ \prod_{\alpha\in V^\star}\left(\frac{f_\alpha}{\lambda_\alpha}\right)^{\lambda_\alpha}.
\]
If, in addition, $V^\star$ is chosen minimally in the sense that the presence of all points in 
$V^\star$ is necessary for $\alpha^\star\in\conv V^\star$ to hold, then 
we also have $\lambda_\alpha>0$ for all $\alpha\in V^\star$. Note that a minimal choice of $V^\star$ is not necessarily unique.

In our setting, due to condition \eqref{eq:c2} the circuit number is a positive number which we associate to each gem degenerate exponent vector of the given polynomial $f$ and, interestingly, besides coercivity questions, the circuit numbers also have been used for analyzing the amoebas, nonnegativity, and sum of squares properties of polynomials supported on circuits (for more details, see \cite{Ili14Wa}).

The derivation of the following necessary conditions for coercivity of a not necessarily gem regular polynomial $f\in\R[x]$ from \cite{BS} bases on a similar technique as presented in \cite{AT,Netz,Rez78}, that is, on evaluations of $f$ along certain curves.

\bthe[Necessary condition for coercivity, \mbox{\cite[Th. 2.29]{BS}}]\label{the:necFirreg}\ \\
Let $f\in\R[x]$ be coercive on $\R^n$. Then the conditions \eqref{eq:c1}-\eqref{eq:c3} are satisfied, and for any $\alpha^\star\in D(f)$ such that
there exists a simplicial face $G\in\mathcal{G}(f)$ with $\alpha^\star\in G$ and $D\cap G=\{\alpha^\star\}$, the following assertions hold.
\i
\item[a)] We have
\q\label{eq:IrregC1}
f_{\alpha^\star}\ \ge\ -\Theta(f,V_G,\alpha^\star).
\qq
\item[b)] For $\alpha^\star\not\in2\N_0^n$ we also have
\q\label{eq:IrregC2}
f_{\alpha^\star}\ \le\ \Theta(f,V_G,\alpha^\star).
\qq
\ii
\ethe

%

Finally we recall the following result from \cite{BS} which, unlike Theorem~\ref{the:char}, guarantees coercivity even for a broad class of gem irregular polynomials.

\bthe[Sufficient condition for coercivity, \mbox{\cite[Th. 3.4]{BS}}]\label{the:suffdeg}\ \\
Let $f\in\R[x]$ be a polynomial satisfying the conditions \eqref{eq:c1}-\eqref{eq:c3}. Furthermore, for each $\alpha^\star\in D(f)$ let $V^\star\subset V(f)$ denote a minimal 
affinely independent set with $\alpha^\star\in\conv V^\star$ and the corresponding unique positive convex coefficients $\{\lambda_{\alpha},\, \alpha\in V^\star\}$ of $\alpha^\star$, let $w(\alpha^\star)>0$, $\alpha^\star\in D(f)$, denote weights with $\sum_{\alpha^\star\in D(f)}w(\alpha^\star)\le1$ and let
\[
f_{\alpha^\star}\ >\ -w(\alpha^\star)\,\Theta(f,V^\star,\alpha^\star)\ \text{ if }\ \alpha^\star\in 2\N_0^n
\]
and 
\[
|f_{\alpha^\star}|\ <\ w(\alpha^\star)\,\Theta(f,V^\star,\alpha^\star)\ \text{ else.}
\] 

Then $f$ is coercive on $\R^n$.
\ethe

\section{Global diffeomorphism property}\label{sec:GlobDiff}

Due to \cite{Gor}, the following theorem, which is of crucial importance for the present work, goes back at least to Jacques S. Hadamard \cite{Had1,Had2,Had3}. For its proof see, e.g., \cite{Gor}, \cite[Sec.~6.2]{Kra}, or \cite[Cor.~4.3]{RP}.

\bthe[Hadamard]\label{thm:Had} A map $F\in C^1(\R^n,\R^n)$ is a $C^1$-diffeomorphism of $\R^n$ onto itself if and only if the map $F$ is proper and $\det JF$ vanishes nowhere on $\R^n$.
\ethe

Since for a continuous map $F:\R^n\rightarrow\R^n$ its properness is equivalent to the property $\|F(x)\|_2^2\rightarrow +\infty$ whenever $\|x\|\rightarrow +\infty$ 
 (see, e.g., \mbox{\cite[Prop. 3.1.15]{Gas}}), one can reformulate Theorem \ref{thm:Had} in the setting of polynomial maps as follows.

\bthe\label{thm:Had2}
A map $F:\R^n\rightarrow\R^n$ with $F=(F_1,\dots,F_n)$, $F_i\in\R[x]$, $i\in I$ is a $C^1$-diffeomorphism of $\R^n$ onto itself if and only if
\q
\det JF(x)\not=0\text{ for all }x\in\R^n\tag{H1}\label{eq:h1}
\qq
and
\q
\text{$\|F(x)\|_2^2\in\R[x]$ is coercive on $\R^n$.}\tag{H2}\label{eq:h2}
\qq
\ethe

In the following we will identify conditions under which \eqref{eq:h1} implies \eqref{eq:h2}, so that the diffeomorphism property of $F$ in Theorem~\ref{thm:Had2} may be characterized by condition \eqref{eq:h1} alone, that is,
the Real Jacobian Conjecture is true under these conditions.
To this end, we shall first show that the function $f:=\|F\|_2^2$ always satisfies the conditions \eqref{eq:c1} and \eqref{eq:c2}.

For any two sets $X_1, X_2 \subseteq\R^n$ we denote by $X_1+X_2:=\{x\in\R^n|\ \exists x_1\in X_1,\,\exists x_2\in X_2: x=x_1+x_2\}\subseteq\R^n$ their Minkowski sum and we define $dX_1:=\{x\in\R^n|\ \exists x_1\in X_1: x=dx_1\}$ for any $d\in\R$. We further denote by $\vertex(P)$ the set of all vertices of some polytope $P\subset\R^n$.
The proof of the following auxiliary result is given in Section \ref{app:lem:Polytopes}.

\begin{lemma}\label{lem:Polytopes} For any polytope $P\subset\R^n$ it holds $v\in \vertex(P+P)$ if and only if $v=2w \text{ with some }w\in\vertex(P).$
\end{lemma}

The subsequent Lemma~\ref{lem:1.1} will provide some useful properties regarding the Newton polytopes at infinity of squared polynomials, while Lemma~\ref{lem:2} shall treat the case of sum of squares polynomials.

\begin{lemma}\label{lem:1.1}

For any $f\in\R[x]$ the following properties hold.
\i
\item[i)] $\new_{\infty}(f^2)=\new_{\infty}(f)+\new_{\infty}(f)$
\item[ii)] $V(f^2)=2V(f)$
\item[iii)] for each $\alpha\in V(f^2)$ it holds $(f^2)_{\alpha}=(f_{\frac{1}{2}\alpha})^2>0$.
\ii
\end{lemma}

\proof
Observe that due to
\[
f^2(x)=\left(\sum_{\alpha\in A(f)}f_{\alpha}x^\alpha\right)^2=\sum_{\alpha,\beta\in A(f)}f_{\alpha}f_{\beta}x^{\alpha+\beta}
\]
\q
=\sum_{\gamma\in A(f)+A(f)}\left(\sum_{\stackrel{\alpha,\beta\in A(f)}{\alpha+\beta=\gamma}}f_\alpha f_\beta\right) x^\gamma,\label{eq:f2}
\qq
the inclusion 
\q
A(f^2)\subseteq A(f)+A(f)\label{eq:inc}
\qq holds, which results in


\[
\new_{\infty}(f^2)=\conv\left(\{0\}\cup A(f^2)\right)\subseteq\conv\left(\{0\}\cup (A(f)+A(f))\right)
\]
\[
\subseteq \conv\left((\{0\}\cup A(f))+(\{0\}\cup A(f))\right)=\conv\left( A_0(f) + A_0(f)\right)
\]
\q
=\conv\left( A_0(f)\right) + \conv\left(A_0(f)\right)=\new_{\infty}(f)+\new_{\infty}(f),\label{eq:inc3}
\qq

where the first inclusion follows from \eqref{eq:inc} and the penultimate equality holds since the convex hull of the Minkowski sum of some given sets is the Minkowski sum of the convex hulls of the sets (see, e.g., \mbox{\cite[Prop. 4.12]{Kra2}}).

Next, we shall show the inclusion 
\q
\new_{\infty}(f)+\new_{\infty}(f)\subseteq \new_{\infty}(f^2)\label{eq:inc4}.
\qq
To this end it suffices to show $\vertex\left(\new_{\infty}(f)+\new_{\infty}(f)\right)\subseteq\new_{\infty}(f^2)$.
By Lemma \ref{lem:Polytopes} we have $\gamma\in\vertex(\new_{\infty}(f)+\new_{\infty}(f))$ if and only if $\gamma=2\delta$ holds with some (unique) $\delta\in\vertex\left(\new_{\infty}(f)\right)=V_0(f)\subseteq A_0(f)$. If $\delta=0$, then $\gamma=0\in\new_{\infty}(f^2)$ by definition. If $\delta\not=0$, then with \eqref{eq:f2} one obtains for the coefficient $(f^2)_{\gamma}\in\R$ of $f^2$ corresponding to the vertex $\gamma$ that
\q
(f^2)_{\gamma}=\sum_{\stackrel{\alpha,\beta\in A(f)}{\alpha+\beta=\gamma}}f_\alpha f_\beta=\left(f_{\delta}\right)^2>0,\label{eq:positive}
\qq
where the last equality holds due to Lemma \ref{lem:Polytopes} and the inequality due to $f_{\delta}\not=0$ following from $0\not=\delta\in A(f)$. This implies $\gamma\in A(f^2)$ and hence $\gamma\in \new_{\infty}(f^2)$. Since $\gamma\in\vertex(\new_{\infty}(f)+\new_{\infty}(f))$ was chosen arbitrarily, the inclusion $\vertex\left(\new_{\infty}(f)+\new_{\infty}(f)\right)\subseteq \new_{\infty}(f^2)$ follows.

The assertion i) follows from \eqref{eq:inc3} and \eqref{eq:inc4}. The assertion ii) follows directly from the assertion i) by using Lemma \ref{lem:Polytopes}. The assertion iii) follows directly from \eqref{eq:positive} above.
\qed

\begin{lemma}\label{lem:2}

For $f(x)=\sum_{i\in I} F_i^2(x)$ with $F_i\in\R[x]$, $i\in I$, the following properties hold.

\i
\item[i)] $\new_{\infty}(f)=\conv\left(\bigcup_{i\in I} 2V_0(F_i)\right)$
\item[ii)] $V_0\left( f \right)\subseteq \bigcup_{i\in I} 2V_0(F_i)$

\item[iii)]$\text{each }\alpha\in V(f)\text{ satisfies }f_{\alpha}>0.$
\ii
\end{lemma}

\proof Observe that due to
\[
f(x)=\sum_{i\in I}F_i^2(x)=\sum_{i\in I}\sum_{\gamma\in A(F_i^2)} (F_i^2)_{\gamma}x^\gamma
\]
\q
=\sum_{\gamma\in\,\bigcup_{i\in I}A(F_i^2)}\,\left(\sum_{i\in I:\,\,\gamma\in A(F_i^2)}(F_i^2)_{\gamma}\right) x^\gamma,\label{eq:sum_1}
\qq
the inclusion
\q
A(f)\subseteq \bigcup_{i\in I} A(F_i^2),\label{eq:inc6}
\qq 
and thus, also the inclusion
\q
A_0(f)\subseteq \bigcup_{i\in I} A_0(F_i^2)\label{eq:inc7}
\qq 
hold.

\textbf{Part i)} 

First, one obtains
\[
\new_{\infty}(f)=\conv\left(A_0(f)\right)\subseteq\conv\left(\bigcup_{i\in I}A_0(F_i^2)\right)= \conv\left(\bigcup_{i\in I}\conv (A_0(F_i^2))\right)
\]
\[
=\conv\left(\bigcup_{i\in I}\new_{\infty}(F_i^2)\right)=\conv\left(\bigcup_{i\in I}\conv (V_0(F_i^2))\right)= \conv\left(\bigcup_{i\in I}V_0(F_i^2)\right) 
\]
\q
=\conv\left(\bigcup_{i\in I}2V_0(F_i)\right),\label{eq:inc8}
\qq
where the inclusion holds due to \eqref{eq:inc7} and the last equality due to Lemma \ref{lem:1.1} ii).
In order to show the other inclusion 
\q
\conv\left(\bigcup_{i\in I}2V_0(F_i)\right)\subseteq\new_{\infty}(f),\label{eq:inclu}
\qq 
it suffices to prove that the vertex set of the polytope $\conv\left(\bigcup_{i\in I}2V_0(F_i)\right)$ is contained in the set $\new_{\infty}(f)$. To this end let $\alpha\in 2\N_0^n$ be an arbitrary vertex of the polytope $\conv\left(\bigcup_{i\in I}2V_0(F_i)\right)$. Then $\alpha$ is necessarily a vertex of each polytope $\conv \left( 2 V_0(F_i)\right)$ containing $\alpha$. Hence, by Lemma \ref{lem:1.1} ii), $\alpha$ is a vertex of each Newton polytope at infinity $\new_{\infty}(F_i^2)$ containing $\alpha$, that is, $\alpha\in V_0(F_i^2)$ for each $i\in I$ with $\alpha\in\new_{\infty}(F_i^2)$. If $\alpha=0$ then obviously $\alpha\in\new_{\infty}(f)$ by definition. Next we shall consider only the case $\alpha\not=0$. Here it holds $\alpha\in V(F_i^2)$ for each $i\in I$ with $\alpha\in\new_{\infty}(F_i^2)$ and using \eqref{eq:sum_1} together with Lemma \ref{lem:1.1} ii) and iii) one obtains
\[
f_{\alpha}=\sum_{i\in I:\,\,\alpha\in A(F_i^2)}\left(F_i^2\right)_{\alpha}=\sum_{i\in I:\,\,\alpha\in V(F_i^2)}\left(F_i^2\right)_{\alpha}
\]
\q
=\sum_{i\in I:\,\,\frac{\alpha}{2}\in V(F_i)}\left((F_i)_{\frac{\alpha}{2}}\right)^2>0.\label{eq:positive2}
\qq
This implies $\alpha\in A(f)$, and hence, $\alpha\in\new_{\infty}(f)$. 


The assertion i) follows from \eqref{eq:inc8} and \eqref{eq:inclu}.

\textbf{Part ii)} 

Due to i) it holds 
\[
V_0(f)=\vertex\left(\text{conv}\left(\bigcup_{i\in I} 2V_0(F_i)\right)\right)\subseteq\bigcup_{i\in I} 2V_0(F_i),
\]
which proves the assertion ii).

\textbf{Part iii)} 

Due to ii) it holds $V(f)\subseteq\bigcup_{i\in I}2V(F_i)$. Thus for any $\alpha\in V(f)$ one has $\alpha\in\bigcup_{i\in I}2V(F_i)$ and with \eqref{eq:positive2} one obtains
\[
f_{\alpha}=\sum_{i\in I:\,\,\frac{\alpha}{2}\in V(F_i)}\left((F_i)_{\frac{\alpha}{2}}\right)^2>0,
\]
which proves the assertion iii).
\qed

The last lemma yields for any sum of squares polynomial $f\in\R[x]$ the following property.
\begin{proposition}\label{prop:C1C2} Every polynomial $f\in\R[x]$ with $f(x)=\sum_{i\in I}F_i^2(x)$, $F_i\in\R[x]$, $i\in I$, fulfills the conditions \eqref{eq:c1} and \eqref{eq:c2}.
\end{proposition}
\proof By Lemma \ref{lem:2} ii) one obtains $V_0\left( f \right)\subseteq \bigcup_{i\in I} 2V_0(F_i)$, which results in
\[
V(f)\subseteq V_0\left( f \right)\subseteq \bigcup_{i\in I} 2V_0(F_i)\subseteq 2\N_0^n,
\]
and thus, $f$ fulfills the condition \eqref{eq:c1}. 

By Lemma \ref{lem:2} iii) one obtains $f_{\alpha}>0$ for each $\alpha\in V(f)$, that is, $f$ also fulfills the condition \eqref{eq:c2}.
\qed

In order to analyze whether the sum of squares polynomial $f=\|F\|_2^2$ corresponding to some polynomial map $F$ also fulfills the condition \eqref{eq:c3}, we shall use the following auxiliary result.

\blem \label{lem:aux} Let $F:\R^n\rightarrow\R^n$ with $F=(F_1,\dots,F_n)$, $F_i\in\R[x]$, $i\in I$, be given. If for each $j\in I$ there exist some $i\in I$ and $k\in\N$ with $ke_{j}\in A(F_i)$, then the polynomial $f=\|F\|_2^2$ satisfies condition \eqref{eq:c3}.
\elem

\proof For every $j\in I$ let there exist some $i\in I$ and some $k\in\N$ such that $ke_{j}\in A(F_{i})$
holds. Define for each $j\in I$ the non-empty set
\[
I(j):=\{i\in I|\ \exists k\in\N\text{ with }ke_j\in A(F_i)\}
\]
and
\[
m(j):=\max\{k\in\N|\ ke_j\in A(F_i),\, i\in I(j)\}
\]
together with the set $\bar{I}(j)\subset I(j)$ of indices at which the maximal value $m(j)$ is attained. For each $j\in I$ it holds $m(j)e_j\in V(F_i)$ for all $i\in \bar{I}(j)$. Using Lemma \ref{lem:1.1} ii) one obtains for each $j\in I$
\[
2m(j)e_j\in V(F_i^2)\quad\text{for all }i\in \bar{I}(j)
\]
and by Lemma \ref{lem:1.1} iii) also
\q\label{eq:ineq-c3}
(F_i^2)_{2m(j)e_j}=\left((F_i)_{m(j)e_j}\right)^2>0.
\qq
With \eqref{eq:sum_1} and \eqref{eq:ineq-c3} one obtains for each $j\in I$
\[
f_{2m(j)e_j}=\sum_{i\in I\,:\, 2m(j)e_j\in A(F_i^2)}(F_i^2)_{2m(j)e_j}=\sum_{i\in\bar{I}(j)}(F_i^2)_{2m(j)e_j}
\]
\[
=\sum_{i\in\bar{I}(j)}\left((F_i)_{m(j)e_j}\right)^2>0,
\] 
which implies $2m(j)e_j\in A(f)$. Since by definition of $m(j)$ it holds $ke_j\notin A(f)$ for all $k>m(j)$, one even obtains that for each $j\in I$ the vector $2m(j)e_j\in A(f)$ is a vertex of $\new_{\infty}(f)$. Thus, we arrive at $2m(j)e_j\in V(f)$ with some $m(j)\in\N$ for every $j\in I$. Thus, $f$ fulfills the condition \eqref{eq:c3}, and the assertion follows.
\qed

In the following, for some vectors $a,b\in\R^n$, we use the notation $a\geq b$ if $a_i\geq b_i$ holds for all $i\in I$, and $\1$ denotes the all-ones vector $(1,\dots,1)^T\in\R^n$.

The next result provides an explicit representation of the Jacobian determinant $\det JF$ of a polynomial map $F$, which will enable us to link the nowhere vanishing property of $\det JF$ to the condition \eqref{eq:c3} of the polynomial $\|F\|_2^2$, as formulated in Proposition \ref{lem:main} below. Its proof is given in Section \ref{app:thm:det}.

\blem[Determinant formula] \label{thm:det}\ \\
Let $F:\R^n\rightarrow\R^n$ with $F=(F_1,\dots,F_n)$, $F_i\in\R[x]$, $i\in I$. Then all $x\in\R^n$ satisfy
\q
\det JF(x)=\sum_{\stackrel{\alpha^i\in A(F_i),\,i\in I}{\sum_{i\in I}\alpha^i\geq \1}}\left(\det(\alpha^1,\ldots,\alpha^n)\,\prod_{i\in I}(F_i)_{\alpha^i}\right)\, x^{\left(\sum_{i\in I}\alpha^i\right) -\1}.\label{det}
\qq
\elem

\bprop\label{lem:main} Let $F:\R^n\rightarrow\R^n$ with $F=(F_1,\dots,F_n)$, $F_i\in\R[x]$, $i\in I$, be given such that 
\[
\det JF(0)\not=0
\] holds. Then the polynomial $f=\|F\|_2^2$ satisfies condition \eqref{eq:c3}.
\eprop

\proof
Assume that $f=\|F\|_2^2$ does not fulfill condition \eqref{eq:c3}. Then by Lemma \ref{lem:aux} there exists an index $j^\star\in I$ such that for every $i\in I$ and every $k\in\N$ one has
$ke_{j^\star}\notin A(F_i)$ 
and, thus, choosing $k=1$ one especially obtains that for all $i\in I$
\q
e_{j^\star}\notin A(F_i)\label{eq:6b}
\qq 
holds. Consider an arbitrary choice of exponent vectors $\alpha^i\in A(F_i)$, $i\in I$, with 
\q\label{eq:suma}
\sum_{i\in I}\alpha^i= \1.
\qq Since $\alpha^i\in\N_0^n$ for each $i\in I$, the system of equations \eqref{eq:suma} implies
\q\label{eq:binary}
\alpha^i_j\in\{0,1\}\text{ for all }i,j\in I.
\qq
Regarding \eqref{eq:suma}, one also has
\q
\sum_{i\in I}\alpha^i_{j^\star}=1
\qq
and thus, due to \eqref{eq:binary}, there exists some (unique) $i^\star\in I$ such that $\alpha^{i^\star}_{j^\star}=1$. By \eqref{eq:6b} there also
exists some $j^{\star\star}\in I\setminus\{j^\star\}$ with $\alpha^{i^\star}_{j^{\star\star}}\not=0$ and, consequently, 
\[
\|\alpha^{i^{\star}}\|_1>1.
\] 
Thus, the binary vector $\alpha^{i^{\star}}$ possesses at least two nonzero entries and, with \eqref{eq:suma}, one obtains
\q
\|\1-\alpha^{i^{\star}}\|_1=\|\sum_{i\in I\setminus\{i^\star\}}\alpha^{i}\|_1<n-1.\label{eq:7}
\qq
By \eqref{eq:7} the remaining $n-1$ binary vectors $\alpha^i$, $i\in I\setminus\{i^\star\}$, can possess at most $n-2$ non-zero entries in total. Thus, by the pigeonhole principle, there exists some $i^{\star\star}\in I\setminus\{i^\star\}$ with
$\alpha^{i^{\star\star}}=0$, 
which results in
\q
\det(\alpha^1,\ldots,\alpha^n)=0.\label{eq:det0}
\qq
Since the choice of vectors $\alpha^i\in A(F_i)$, $i\in I$, with \eqref{eq:suma} was arbitrary, using Lemma \ref{thm:det} and \eqref{eq:det0} one finally obtains
\[
\det JF(0)=\sum_{\substack{\alpha^i\in A(F_i),\,i\in I\\\sum_{i\in I}\alpha^i= \1}}\left(\det(\alpha^1,\ldots,\alpha^n)\,\prod_{i=1}^n(F_i)_{\alpha^i}\right)=0,
\]
and the assertion follows.

\qed

The combination of Propositions \ref{prop:C1C2} and  \ref{lem:main} provides the following result.

\begin{proposition}\label{prop:main} Let $F:\R^n\rightarrow\R^n$ with $F=(F_1,\dots,F_n)$, $F_i\in\R[x]$, $i\in I$, be given such that 
\[
\det JF(0)\not=0
\] 
holds. 
Then the polynomial $f=\|F\|_2^2$ fulfills the conditions \eqref{eq:c1}-\eqref{eq:c3}.
\end{proposition}

The following two theorems contain the main results of this article. The first one assumes the gem-regularity of the polynomial $\|F\|_2^2$, while the second one treats also the case of gem-irregular polynomials $\|F\|_2^2$ under some further conditions imposed on the coefficients corresponding to the gem-degenerate exponent vectors of $\|F\|_2^2$ which also include the circuit number information.

\bthe\label{thm:main} For $F:\R^n\rightarrow\R^n$ with $F=(F_1,\dots,F_n)$, $F_i\in\R[x]$, $i\in I$, let the polynomial $f=\|F\|_2^2$ be gem regular. Then the following two assertions are equivalent. 
\i
\item[a)] $F$ is a $C^1$-diffeomorphism of $\R^n$ onto itself.
\item[b)] $\det JF(x)\not=0 \text{ holds for all }x\in\R^n$.
\ii
\ethe

\proof Assertion a) implies b) by direct application of Theorem \ref{thm:Had}. For the proof of the reverse direction observe that assertion b) and Proposition \ref{prop:main} imply that $\|F\|_2^2$ fulfills the conditions \eqref{eq:c1}-\eqref{eq:c3} which, by Theorem \ref{the:char}, characterize the coercivity on $\R^n$ of the gem regular polynomial $\|F\|_2^2$. The map $F$ thus fulfills the conditions  \eqref{eq:h1} and \eqref{eq:h2}, and Theorem~\ref{thm:Had2} finally implies that $F$ is a $C^1$-diffeomorphism of $\R^n$ onto itself.
\qed

%

\bthe\label{thm:maindeg} For $F:\R^n\rightarrow\R^n$ with $F=(F_1,\dots,F_n)$, $F_i\in\R[x]$, $i\in I$, let $f=\|F\|_2^2$. For each $\alpha^\star\in D(f)$ let $V^\star\subset V(f)$ denote a minimal affinely independent set with $\alpha^\star\in\conv V^\star$ and the corresponding unique positive convex coefficients $\lambda_{\alpha}$, $\alpha\in V^\star$, of $\alpha^\star$, let $w(\alpha^\star)>0$, $\alpha^\star\in D(f)$, denote weights with $\sum_{\alpha^\star\in D(f)}w(\alpha^\star)\le1$, and let further
\[
f_{\alpha^\star}\ >\ -w(\alpha^\star)\,\Theta(f,V^\star,\alpha^\star)\ \text{ if }\ \alpha^\star\in 2\N_0^n
\]
as well as
\[
|f_{\alpha^\star}|\ <\ w(\alpha^\star)\,\Theta(f,V^\star,\alpha^\star)\ \text{ else.}
\] 
Then the following two assertions are equivalent.
\i
\item[a)] $F$ is a $C^1$-diffeomorphism of $\R^n$ onto itself.
\item[b)] $\det JF(x)\not=0 \text{ holds for all }x\in\R^n$.
\ii
\ethe

\proof The proof runs along the same lines as the proof of Theorem \ref{thm:main}, where Theorem \ref{the:char} is replaced by Theorem \ref{the:suffdeg}.
\qed

\begin{example}\label{ex:thm2} Consider the polynomial map $F_t:\R^2\rightarrow\R^2$ with $F_{t,1}(x)=x_1+x_1^3-tx_2^3$ and $F_{t,2}(x)=x_2+x_1^3+x_2^3$ for some parameter value $t\in\R$. We shall show that the map $F_t$ is a $C^1$-diffeomorphism of $\R^2$ onto itself for all parameter values $t>-1$, and that $F_t$ does not possess this diffeomorphism property for any $t<-1$. 

We define $f_t(x):=\|F_t(x)\|_2^2$ for all $x\in\R^2$. First, let us consider the case $t=1$. Observe that 
\[
\det JF_1(x)=1+3x_1^2+3x_2^2+18x_1^2x_2^2>0\quad\text{ for all $x\in\R^2$}
\] 
holds and hence $\det JF_1(x)\not=0$ for all $x\in\R^2$ is fulfilled. One further obtains
\[
f_1(x)=2x_1^6+2x_2^6+2x_1^4+2x_2^4+2x_1^3x_2-2x_1x_2^3+x_1^2+x_2^2
\]
with the corresponding gem
\[
\G(f_1)=\conv\left((6,0),(0,6)\right),
\]
which implies gem regularity of $f_1$, since $D(f_1)=V^c(f_1)\cap\G(f_1)=\emptyset$ (see Fig. \ref{Fig:F2}). According to Theorem \ref{thm:main} the map $F_1$ thus is a $C^1$-diffeomorphism of $\R^2$ onto itself.

Next we shall consider only parameter values $t\not=1$. First, observe that the condition
\q
\det JF_t(x)=1+3x_1^2+3x_2^2+(9+9t)x_1^2x_2^2\not=0 \text{ for all }x\in\R^2\label{eq:notvanish}
\qq
is violated for any $t<-1$, since the choice $x(s)=(s,s)$ with $s\in\R$ leads to the function $\det JF_t(x(s))=1+6s^2+(9+9t)s^4$ which 
possesses real zeros. By Theorem~\ref{thm:Had}, $F_t$ can thus not be a $C^1$-diffeomorphism of $\R^2$ onto itself for any $t<-1$.

On the other hand, \eqref{eq:notvanish} holds for all $t\geq -1$, since the Jacobian determinant then is strictly positive. One further obtains
\[
f_t(x)=2x_1^6+(1+t^2)x_2^6+2(1-t)x_1^3x_2^3+2(x_1^4+x_2^4)+2(x_1^3x_2-tx_1x_2^3)+x_1^2+x_2^2
\]
with
\[
\G(f_t)=\conv\left((6,0),(0,6)\right),
\]
which, for parameter values $t\not=1$, implies gem irregularity of $f_t$, since $D(f_t)=V^c(f)\cap\G(f_t)=\{(3,3)\}$ due to $f_{t,{(3,3)}}=2(1-t)\not=0$ (see Fig.~\ref{Fig:F2}). For $\alpha^\star=(3,3)\in D(f_t)$ the unique minimal affinely independent subset $V^\star\subseteq V(f_t)$ with $\alpha^\star\in \conv V^\star$ is given by the vertex set at infinity of $\new_{\infty}(f_t)$ itself, that is, $V^\star:=V(f_t)=\{(6,0),(0,6)\}$. From the convex representation $\alpha^\star=(3,3)=\frac{1}{2}(6,0)+\frac{1}{2}(0,6)$  with the unique positive convex coefficients $\lambda_{(6,0)}=\lambda_{(0,6)}=\frac{1}{2}$, computing the corresponding circuit number yields
\[
\Theta(f_t,\alpha^\star,V(f_t))=\left(\frac{f_{t,{(6,0)}}}{\lambda_{(6,0)}}\right)^{\lambda_{(6,0)}}\left(\frac{f_{t,{(0,6)}}}{\lambda_{(0,6)}}\right)^{\lambda_{(0,6)}}=\sqrt{\frac{2}{\frac{1}{2}}}\sqrt{\frac{1+t^2}{\frac{1}{2}}}
\] 
\q\label{eq:circuit}
=2\sqrt{2}\sqrt{(1+t^2)}.
\qq
Further, choosing the weight $w((3,3)):=1$, the inequality
\[
|f_{t,{(3,3)}}|<\Theta(f_t,\alpha^\star,V(f_t))
\]
holds if and only if $t\not=-1$, because due to $f_{t,{(3,3)}}=2(1-t)$ and \eqref{eq:circuit}, the inequality
\[
|2(1-t)|<2\sqrt{2}\sqrt{(1+t^2)}
\]
holds if and only if $t\not=-1$. According to Theorem \ref{thm:maindeg}, the latter fact together with  \eqref{eq:notvanish} imply that the map $F_t$ is a $C^1$-diffeomorphism of $\R^2$ onto itself for all parameter values $t>-1$, and the assertion follows.

\begin{figure}[htbp]
\begin{center}
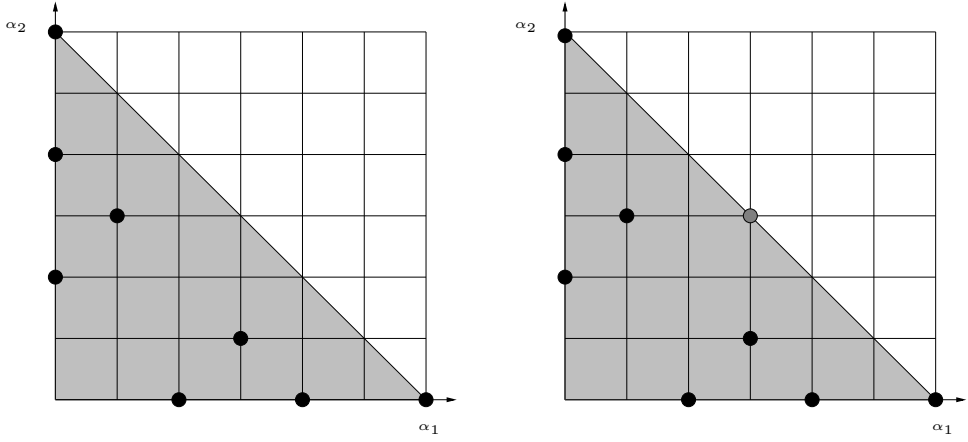
\caption{Illustration of Example \ref{ex:thm2}. In the left picture, the shaded area corresponds to the Newton polytope at infinity $\new_{\infty}(f_1)$, and the black circles stand for the set $A(f_1)$. For the case $t\not=1$, in the right picture, the shaded area corresponds to the Newton polytope at infinity $\new_{\infty}(f_t)$, the black circles stand for the set $A(f_t)\setminus D(f_t)$, and the shaded circle describes the (singleton) set $D(f_t)$.}
\label{Fig:F2}
\end{center}
\end{figure}


\end{example}

\begin{remark}
Example~\ref{ex:thm2} also shows that, despite the assertion of Proposition~\ref{lem:main}, under the assumptions of Theorems~\ref{thm:main} or \ref{thm:maindeg} the diffeomorphism property of a polynomial map $F$ may not solely be 
characterized by the condition $\det JF(0)\neq0$. In fact, in the example we have $\det JF_t(0)\neq0$ for any $t\in\R$, but $F_t$ is not a $C^1$-diffeomorphism
for $t<-1$.
\end{remark}

\begin{remark}
In \cite[Lem. 2.22]{BS} we showed that gem regularity of a polynomial $f$ is a weak condition in the sense that it follows from a general position property of the multiplier vectors $\alpha\in A(f)$. Unfortunately,
the polynomials $f=\|F\|_2^2$ considered in the present paper possess a special structure, so that gem regularity of such functions is not necessarily a mild assumption. 

In fact, Example~\ref{ex:thm2} provides a parametric
family of such polynomials for which gem regularity and Theorem~\ref{thm:main} may only be employed at a single choice of the parameter ($t=1$). On the other hand, Theorem~\ref{thm:maindeg} covers the gem irregular case well enough to treat all members of the
parametric family except for a singular choice of the parameter ($t=-1$), which will be considered separately in Example~\ref{ex:thm2_voll} below.
\end{remark}

\section{Coercivity under linear transformations}\label{sec:lin}

In this section we shall show how linear transformations can help to study the global diffeomorphism property of a polynomial map  when the assumptions of Theorems~\ref{thm:main} and \ref{thm:maindeg} are violated, like for the singular parameter value in Example~\ref{ex:thm2}.

\begin{proposition}\label{thm:Coerc_LCT}
For any regular matrix $A\in\R^{n\times n}$ a function $f:\R^n\to\R$ is coercive on $\R^n$ if and only if the function $f\circ A^{-1}$ is 
coercive on $\R^n$.
\end{proposition}

\proof  Let $A$ be a regular matrix, let $f$ be coercive on $\R^n$, and consider a sequence $(y^\nu)\subset\R^n$ with $\lim_{\nu\to\infty}\|y^\nu\|=+\infty$. Then 
we have
$\|y^\nu\|\leq \|A\|\|A^{-1}y^\nu\|$ for all $\nu\in\N$, where $\|A\|$ denotes the matrix norm of $A$ induced by $\|\cdot\|$.
This implies $\lim_{\nu\to\infty}\|A^{-1}y^\nu\|=+\infty$ and, by the coercivity of $f$, $\lim_{\nu\to\infty}f(A^{-1}y^\nu)=+\infty$, 
so that the coercivity of $f\circ A^{-1}$  is shown. The reverse direction may be shown along the same lines, using the 
identity $f=(f\circ A^{-1})\circ A$.
\qed

We may also improve the formulation of Proposition~\ref{prop:main} as follows:
\begin{proposition}\label{prop:mainA} Let $F:\R^n\rightarrow\R^n$ with $F=(F_1,\dots,F_n)$, $F_i\in\R[x]$, $i\in I$, be given such that 
\[
\det JF(0)\not=0
\] 
holds. 
Then for any regular matrix $A\in\R^{n\times n}$ the polynomial $\|F\circ A^{-1}\|_2^2$ fulfills the conditions \eqref{eq:c1}-\eqref{eq:c3}.
\end{proposition}
\proof Since $F\circ A^{-1}$ is a polynomial map, the assertion follows from Proposition \ref{prop:main} and
$
\det J(F\circ A^{-1})(0)=\det JF(0)\cdot\det A^{-1}\neq0.
$
\qed

\bcor\label{cor:main}
The assertion of Theorem~\ref{thm:main} remains true, if the assumption of gem regularity of the polynomial $\|F\|_2^2$ is replaced by the assumption of gem regularity of the polynomial $\|F\circ A^{-1}\|_2^2$ for some regular matrix $A\in\R^{n\times n}$.
\ecor
\proof We only have to modify the proof that assertion b) implies assertion a).
In fact, for the given matrix $A\in\R^{n\times n}$ assertion b) and Proposition \ref{prop:mainA} imply that $\|F\circ A^{-1}\|_2^2$ fulfills the conditions \eqref{eq:c1}-\eqref{eq:c3} which, by Theorem \ref{the:char}, characterize the coercivity on $\R^n$ of the gem regular polynomial $\|F\circ A^{-1}\|_2^2$. Consequently, by Proposition~\ref{thm:Coerc_LCT} also the polynomial $\|F\|_2^2$ is coercive, so that the map $F$ thus fulfills the conditions  \eqref{eq:h1} and \eqref{eq:h2}, and Theorem~\ref{thm:Had2} implies the assertion.
\qed

The following result is shown analogously.

\bcor\label{cor:maindeg}
The assertion of Theorem~\ref{thm:maindeg} remains true, if the assumptions on the polynomial $\|F\|_2^2$ are replaced by the same assumptions on the polynomial $\|F\circ A^{-1}\|_2^2$ for some regular matrix $A\in\R^{n\times n}$.
\ecor


We illustrate Corollary~\ref{cor:maindeg} by sharpening the result given in Example \ref{ex:thm2}.

\begin{example}\label{ex:thm2_voll} Consider the polynomial map $F_t:\R^2\rightarrow\R^2$ from Example~\ref{ex:thm2} with $F_{t,1}(x)=x_1+x_1^3-tx_2^3$ and $F_{t,2}(x)=x_2+x_1^3+x_2^3$ for some parameter value $t\in\R$. Then $F_t$ is a $C^1$-diffeomorphism of $\R^2$ onto itself if and only if $t\geq -1$. 

In fact, by Example \ref{ex:thm2}, it suffices to show that $F_t$ is a $C^1$-diffeomorphism of $\R^2$ onto itself for the singular parameter value $t=-1$. In fact, using the linear coordinate transformation $x=A^{-1}y$ with the matrix
\[
A^{-1}=\begin{pmatrix} 1 & 1\\ 1 & -1 \end{pmatrix},
\]
one obtains the gem-irregular polynomial
\qa
f_{-1}(A^{-1}y) & := &\|F_{-1}(A^{-1}y)\|_2^2=
\qqa
\qa
=2y_1^2+2y_2^2+8y_1^4+24y_1^2y_2^2+8y_1^6+48y_1^4y_2^2+72y_1^2y_2^4
\qqa
with $D(f_{-1}\circ A^{-1})=\{(4,2)\}$ and $V(f_{-1}\circ A^{-1})=\{(6,0),(2,4),(0,2)\}$. From the positivity of the circuit number $\Theta(f_{-1}\circ A^{-1},(4,2),V(f_{-1}\circ A^{-1}))$ corresponding to the unique gem-degenerate exponent vector $\alpha^\star=(4,2)$ of $f_{-1}\circ A^{-1}$, one obtains the inequality
\qa\label{eq:circuit3}
48=\left(f_{-1}\circ A^{-1}\right)_{(4,2)}>-\Theta(f_{-1}\circ A^{-1},(4,2),V(f_{-1}\circ A^{-1})).
\qqa
Since $\det JF_{-1}(x_1,x_2)=1+3x_1^2+3x_2^2>0$ holds for all $x\in\R^2$, Corollary~\ref{cor:maindeg} yields the assertion.
\end{example}


\begin{remark}
In Pinchuk's counterexample to the Real Jacobian Conjecture (cf. \cite{Pi}), the Jacobian determinant of $F$ vanishes nowhere on $\R^2$ so that, by Proposition~\ref{prop:mainA}, the sum of squares polynomial
$\|F\circ A^{-1}\|_2^2$ does satisfy the conditions \eqref{eq:c1}-\eqref{eq:c3}
for any regular matrix $A\in\R^{2\times 2}$. Since, however, $F$ is not a global $C^1$-diffeomorphism, $\|F\circ A^{-1}\|_2^2$ can neither be gem regular nor satisfy the additional sufficient conditions from 
Theorem~\ref{thm:maindeg} for any regular matrix $A\in\R^{2\times 2}$.
\end{remark}

\section{Final remarks}\label{sec:fin}

This article shows that the global diffeomorphism property of a real polynomial map $F:\R^n\to\R^n$ can sometimes be studied by analyzing the coercivity property of the sum of squares polynomial $\|F\|_2^2$ via its Newton polytope at infinity $\new_{\infty}(\|F\|_2^2)$. However, due to the special structure of the polynomial $\|F\|_2^2$, the assumptions of known sufficiency theorems for coercivity are not necessarily mild and may be expected to be violated. On the other hand, while preserving the coercivity property, suitable linear coordinate transformations may help to transform such a degenerated polynomial into another one, for which the known techniques for verifying coercivity can be applied. 

In order to better understand the coercivity property of multivariate polynomials over $\R^n$, it is thus an interesting question whether for each coercive polynomial $f$ there exists some linear coordinate transformation such that, in new coordinates, $f$ fulfills the conditions from Theorem~\ref{the:char} or from Theorem~\ref{the:suffdeg}, and how such a linear coordinate transformation may be constructed. We leave these questions for future research.

\section*{Acknowledgments} The authors are grateful to Yu. Nesterov and V. Shikhman for pointing out the importance of the invariance of coercivity under linear transformations, and for other fruitful discussions on the subject of this article.


\renewcommand{\thesection}{A}
\section{Appendix}
\setcounter{definition}{0}

\subsection{Proof of Lemma \ref{lem:Polytopes}}\label{app:lem:Polytopes}

For any $\bar{v}\in\vertex(P+P)$ there exists a vector $a\in\R^n\setminus\{0\}$ such that $\bar{v}$ is the unique optimal point of the problem
\[
\text{(LP1)}\quad\max_{v\in P+P}a^Tv.
\]
Let $\bar{w}\in \vertex(P)$ be an optimal point of the problem
\[
\text{(LP2)}\quad\max_{w\in P}a^Tw.
\]
Since 
\[
\max_{v\in P+P}a^Tv=\max_{(x,y)\in P\times P}a^T(x+y)=\max_{x\in P}a^Tx+\max_{y\in P}a^Ty
\]
\q
=2\max_{x\in P}a^Tx=2a^T\bar{w}=a^T2\bar{w}\label{eq:Polytope1}
\qq
holds, the point $2\bar{w}\in 2\vertex(P)$ is an optimal point of (LP1). A $a$ was chosen such that $\bar{v}$ is  the unique optimal point of (LP1), one obtains $\bar{v}=2\bar{w}$ with $\bar{w}\in\vertex(P)$.


On the other hand, choose $\bar{w}\in \vertex(P)$ and put $\bar{v}=2\bar{w}$. To show is $\bar{v}\in \vertex(P+P)$. Observe that there exists some $a\in\R^n\setminus\{0\}$ such that $\bar{w}$ is the unique optimal point of the problem (LP2). Using \eqref{eq:Polytope1}, the point $\bar{v}$ is thus an optimal point of (LP1). Assume that $\bar{v}\notin\vertex(P+P)$ holds. Since (LP1) must possess a vertex solution, there exists an optimal point $\bar{z}:=\bar{x}+\bar{y}\in P+P$ of (LP1) with $\bar{v}\not=\bar{z}$. For the point $\bar{u}:=\textstyle{\frac12}(\bar{x}+\bar{y})\in P$ we obtain the identity
\[
a^T\bar{u}=
\frac{1}{2}a^T\left(\bar{x}+\bar{y}\right)=\frac{1}{2}a^T\bar{z}=a^T\bar{w},
\]
where the last equation holds since both $\bar{z}$ and $\bar{v}=2 \bar{w}$ are optimal for (LP1).
The point $\bar{u}\in P$ is thus an optimal point of the problem (LP2), and the uniqueness of $\bar{w}$ implies $\bar{w}=\bar{u}$. This leads to the contradiction $\bar{v}=\bar{z}$, and thus the assertion $\bar{v}\in \vertex(P+P)$ follows.
\qed


\subsection{Proof of Lemma \ref{thm:det}}\label{app:thm:det}

Let $S_n$ denote the symmetric group on $n$ elements, let $\sign(\sigma)$ denote the permutation sign of $\sigma\in S_n$,
and for some arbitrarily given $x\in\R^n$ let the entries of $JF(x)$ be denoted by $a_{ij}$, $i,j\in I$.
Then the Leibniz formula for determinants yields
\[
\det JF(x)\ =\ \sum_{\sigma\in S_n}\sign(\sigma)\prod_{i\in I} a_{i,\sigma(i)}
\]
with
\[
a_{i,\sigma(i)}\ =\ \frac{\partial}{\partial x_{\sigma(i)}}F_i(x)\ =\ \sum_{\alpha^i\in A(F_i)}(F_i)_{\alpha^i}\,\frac{\partial}{\partial x_{\sigma(i)}}x^{\alpha^i}
\]
for all $\sigma\in S_n$ and $i\in I$. Interchanging multiplication and addition, and splitting the appearing products, further leads to
\[
\prod_{i\in I} a_{i,\sigma(i)}\ =\ \sum_{(\alpha^1,\dots,\alpha^n)\in A(F_1)\times\dots\times A(F_n)}\left[\prod_{i\in I}(F_i)_{\alpha^i}\right]\cdot\left[\prod_{i\in I}\frac{\partial}{\partial x_{\sigma(i)}}x^{\alpha^i}\right]
\]
for all $\sigma\in S_n$. In fact, in the above summation for any $i\in I$ it is sufficient to choose $\alpha^i\in A(F_i)$ with $\alpha^i_{\sigma(i)}\ge1$,
since the existence of some $j\in I$ with $\alpha^j_{\sigma(j)}=0$ means that the monomial $x^{\alpha^j}$ does not depend on the variable $x_{\sigma(j)}$, resulting in
\[
\frac{\partial}{\partial x_{\sigma(j)}}x^{\alpha^j}\ =\ 0\quad\text{and}\quad \prod_{i\in I}\frac{\partial}{\partial x_{\sigma(i)}}x^{\alpha^i}=0.
\]
This shows
\[
\prod_{i\in I} a_{i,\sigma(i)}\ =\ \sum_{\substack{\alpha^i\in A(F_i),\,i\in I\\\alpha^i_{\sigma(i)}\ge1,\,i\in I}}\left[\prod_{i\in I}(F_i)_{\alpha^i}\right]\cdot\left[\prod_{i\in I}\frac{\partial}{\partial x_{\sigma(i)}}x^{\alpha^i}\right]
\]
for all $\sigma\in S_n$. 

Next, for any $(\alpha^1,\ldots,\alpha^n)$ in the above summation and any $i\in I$ we have
\[
\frac{\partial}{\partial x_{\sigma(i)}}x^{\alpha^i}\ =\ \alpha^i_{\sigma(i)}x^{\alpha^i_{\sigma(i)}-1}_{\sigma_i}\prod_{j\neq i}x^{\alpha^j_{\sigma(j)}}_{\sigma(j)}
\]
and, since $\sigma$ is a permutation,
\[
\prod_{i\in I}\frac{\partial}{\partial x_{\sigma(i)}}x^{\alpha^i}\ =\ \left[\prod_{i\in I}\alpha^i_{\sigma(i)}\right]\cdot x^{\sum_{i\in I}\alpha^i-\1}.
\]
We arrive at
\qa
\prod_{i\in I} a_{i,\sigma(i)} & = & \sum_{\substack{\alpha^i\in A(F_i),\,i\in I\\\alpha^i_{\sigma(i)}\ge1,\,i\in I}}\left[\prod_{i\in I}(F_i)_{\alpha^i}\right]\cdot\left[\prod_{i\in I}\alpha^i_{\sigma(i)}\right]\cdot x^{\sum_{i\in I}\alpha^i-\1}\\
& = & \sum_{\substack{\alpha^i\in A(F_i),\,i\in I\\\alpha^i_{\sigma(i)}\ge1,\,i\in I}}\left[\prod_{i\in I}\alpha^i_{\sigma(i)}\right]\cdot m(\alpha^1,\ldots,\alpha^n,x)
\qqa
for all $\sigma\in S_n$, where the monomial
\[
m(\alpha^1,\ldots,\alpha^n,x)\ :=\ \left[\prod_{i\in I}(F_i)_{\alpha^i}\right]\cdot x^{\sum_{i\in I}\alpha^i-\1}
\]
does not depend on $\sigma$. Hence, we may write
\qa
\det JF(x) & = & \sum_{\sigma\in S_n}\sign(\sigma)\sum_{\substack{\alpha^i\in A(F_i),\,i\in I\\\alpha^i_{\sigma(i)}\ge1,\,i\in I}}\left[\prod_{i\in I}\alpha^i_{\sigma(i)}\right]\cdot m(\alpha^1,\ldots,\alpha^n,x)\\
& = & \sum_{\substack{\alpha^i\in A(F_i),\,i\in I}}\sum_{\substack{\sigma\in S_n\\\alpha^i_{\sigma(i)}\ge1,\,i\in I}}\sign(\sigma)\cdot\left[\prod_{i\in I}\alpha^i_{\sigma(i)}\right]\cdot m(\alpha^1,\ldots,\alpha^n,x).
\qqa
In the latter outer summation it suffices to consider $\alpha^i\in A(F_i)$, $i\in I$, with $\sum_{i\in I}\alpha^i\ge\1$, since otherwise there would exist some $j\in I$ with $\alpha^i_j=0$ for all $i\in I$,
resulting in $\alpha^{\sigma^{-1}(j)}_j=0$ for any $\sigma\in S_n$. However, then the inner summation would be taken over the empty set.

After introducing this restriction on the outer summation, we may drop the constraint $\alpha^i_{\sigma(i)}\ge1$, $i\in I$, in the inner summation since, for given $\sigma\in S_n$, its violation leads to
a vanishing product $\prod_{i\in I}\alpha^i_{\sigma(i)}$. Thus we have shown the assertion
\qa
\det JF(x) & = & \sum_{\substack{\alpha^i\in A(F_i),\,i\in I\\\sum_{i\in I}\alpha^i\ge\1}}\sum_{\sigma\in S_n}\sign(\sigma)\cdot\left[\prod_{i\in I}\alpha^i_{\sigma(i)}\right]\cdot m(\alpha^1,\ldots,\alpha^n,x)\\
& = & \sum_{\substack{\alpha^i\in A(F_i),\,i\in I\\\sum_{i\in I}\alpha^i\ge\1}}\det(\alpha^1,\ldots,\alpha^n)\cdot m(\alpha^1,\ldots,\alpha^n,x),
\qqa
where the final identity is due to the Leibniz formula for determinants.\qed

\end{document}